\documentclass[twocolumn,10pt]{article}

\usepackage{amsmath,amssymb,amsfonts}
\usepackage{algorithmic}
\usepackage{graphicx}
\usepackage{textcomp}
\usepackage[dvipsnames,svgnames,x11names]{xcolor}
\usepackage{mathptmx}
\usepackage[scaled=.81]{helvet}
\usepackage{courier}
\usepackage{newpxmath}
\usepackage[bf,compact,small]{titlesec}
\titleformat*{\section}{\sffamily\bfseries\color{NavyBlue}}
\titleformat*{\subsection}{\sffamily\bfseries}
\titleformat*{\subsubsection}{\sffamily\bfseries}
\titleformat*{\paragraph}{\sffamily\bfseries}
\usepackage[margin=9ex]{geometry}

\newcommand{\runninghead}[1]{}
\newcommand{\affiliation}[1]{\footnotetext[0]{\sf #1}}
\def\affilnum#1{${}^{\text{{#1}}}$}
\newcommand{\keywords}[1]{\textbf{Keywords} #1}
\newenvironment{acks}{\section*{Acknowledgements}}{}
\newcommand{\citep}[1]{\cite{#1}}
\newcommand{\corrauth}[1]{}
\newcommand{\email}[1]{}

\begin{document}

\newcommand{\blasplusplus}{BLAS\texttt{++}\xspace}
\newcommand{\cplusplus}{C$\!$\texttt{++}\xspace}
\newcommand{\lapackplusplus}{LAPACK$\!$\texttt{++}\xspace}
\newcommand{\openmp}{\textsf{OpenMP}\xspace}

\runninghead{Portable Performance for Numerical}

\title{HPL-MxP Benchmark: Mixed-Precision Algorithms, Iterative
Refinement, and Scalable Data Generation}

\author{Jack Dongarra\affilnum{1, 2, 3} %
and Piotr Luszczek\affilnum{4, 1}}

\date{September 23, 2025\thanks{Use of this work is controlled by the %
human-to-human license listed in Exhibit 3 of %
https://doi.org/10.48550/arXiv.2306.09267}}

\maketitle

\affiliation{%
\affilnum{1} University of Tennessee, %
\affilnum{2} Oak Ridge National Laboratory, %
\affilnum{3} University of Manchester, %
\affilnum{4} MIT Lincoln Laboratory} %
\corrauth{Piotr Luszczek}
\email{piotr.luszczek@ll.mit.edu}

\begin{abstract}
We present a mixed-precision benchmark called HPL-MxP that uses both a
lower-precision LU factorization with a non-stationary iterative
refinement based on GMRES. We evaluate the numerical stability of one of
the methods of generating the input matrix in a scalable fashion and
show how the diagonal scaling affects the solution quality in terms of
the backward-error. Some of the performance results at large scale
supercomputing installations produced Exascale-level compute throughput
numbers thus proving the viability of the proposed benchmark for
evaluating such machines. We also present the potential of the benchmark
to continue increasing its use with proliferation of hardware
accelerators for AI workloads whose reliable evaluation continues to
pose a particular challenge for the users.
\end{abstract}

\keywords{%
floating point representation, %
hardware accelerator offloading, %
high performance computing, %
numerical linear algebra, %
mixed-precision algorithms %
}

\begin{acks}
Research was sponsored by the Department of the Air Force Artificial
Intelligence Accelerator and was accomplished under Cooperative
Agreement Number FA8750-19-2-1000. The views and conclusions contained
in this document are those of the authors and should not be interpreted
as representing the official policies, either expressed or implied, of
the Department of the Air Force or the U.S.\ Government. The U.S.\
Government is authorized to reproduce and distribute reprints for
Government purposes notwithstanding any copyright notation herein.

This research uses resources of the Oak Ridge Leadership Computing
Facility, which is a DOE Office of Science User Facility supported under
Contract DE-AC05-00OR22725. This research also uses resources of the
Argonne Leadership Computing Facility, which is a DOE Office of
Science User Facility supported under Contract DE-AC02-06CH11357.

This work was partially funded by an NSF Collaborative Research
Framework: Basic ALgebra LIbraries for Sustainable Technology with
Interdisciplinary Collaboration (BALLISTIC), a project of the
University of Tennessee 2004541.
\end{acks}


\section{Introduction}
\label{sec:intro}

The use of mixed-precision, while now prevalent across scientific and
artificial intelligence (AI) applications, is not necessarily uniform in
terms of ubiquity nor in terms of relative performance gains. As we show
in the references, theory, and results below, there are gaps dependent
on the discipline and data properties with the larger scientific
application ecosystem.
Updating these gaps in terms of simulations'
needs relies on the fundamental software
libraries that started offering some of the mixed-precision
functionality: some of the solvers are available
with proven theoretical bounds on
the quality of the results.  It is important
to point out that the protracted
sunset of Moore's Law~\citep{moore1965cramming} created a new avenue in the pursuit
of the improvements in hardware performance, namely limiting the number
of bits of the processed data thus leading to not only new
floating-point precisions but also much skewed throughput rates for
low-bit formats, either fixed- or floating-point. This started creating
performance inversions for higher precision data as is the case for the
upgrade from NVIDIA Hopper to Blackwell GPUs that lowers the performance
of the IEEE standard 64-bit format.

In this environment, the HPL-MxP benchmark (formerly known as HPL-AI)
seeks to highlight these emerging convergent trends and provide a
unified framework for testing hardware capabilities relevant to both
high-performance computing (HPC) and AI workloads. On the one hand, the
traditional HPC continues focusing on high-fidelity simulations for
modeling phenomena in physics, chemistry, biology, and so on. The
specialized mathematical models that drive these computational endeavors
require, for the majority of their data sets, 64-bit accuracy. On the
other hand, the machine learning and data science methods that fuel
advances in AI manage to achieve the desired results at 32-bit and now
more often in even lower floating-point precision formats. This
drastically lesser demand for computational accuracy fueled a resurgence
of interest in both new hardware platforms and software implementations
of new numerical methods. In combination, both of
these advances and the
HPL-MxP benchmark aim to showcase a mix of unprecedented performance
levels and accompanying energy savings to achieve the simulation,
classification, and recognition fidelity thought to be only possible
with higher-accuracy data formats.

To resolve this newly emergent directions in hardware and software
development, HPL-MxP aims to strike a balance between the divergent
trends by delivering a combination of modern algorithms and efficient
use of modern accelerator hardware while simultaneously drawing a
connection to the established solver inside the decades-old
HPL benchmark~\citep{Dongarra2003} and its deployment base of the largest
supercomputing installations in the world. The method of choice used in
the primary solver is a judicious use of the LU factorization in lower
precision and the iterative refinement performed in a manner that brings
the solution back to the 64-bit accuracy expected from the HPL results.
The main innovation of HPL-MxP (in contrast to HPL) lies in dropping the
requirement of 64-bit computation throughout the entire computation and
instead opting for low-precision (likely 16-bit or lower) accuracy
inside LU, and a sophisticated refinement to recover the accuracy lost
during the factorization. The iterative method guaranteed to be
numerically stable is the generalized minimal residual method (GMRES),
which uses application of the L and U factors to serve as a
preconditioner. The combination of these algorithms is demonstrably
sufficient for high accuracy~\citep{carson2018iterref,carson2017iterref}
and may be implemented in a way that takes advantage of the current and
upcoming devices for accelerating AI
workloads~\citep{haidar2018mxpiterref,haidar2020mxpiterref}.



\section{Background and Historical Perspective}
\label{sec:bckgrnd}

Backward numerical stability is a method developed in
the 1960s~\citep{wilkinson1963rea,wilkinson1965aep}
and together with iterative
refinement in higher precision was to be used to improve the solution to
a system of linear equations $Ax=b$ beyond the accuracy achieved in a
uniform precision scenario. Many methods, either direct or iterative,
are able to achieve low values for $\|Ax-b\|$ such that it can be
bound by the corresponding norms of $A$, $x$, and $b$ scaled by a
slow-growing function of $n$ or the size of the matrix $A$ as well as
the vectors $x$ and $b$.
In this classical view, the iterative refinement was only an optional
component because the original solution was sufficiently accurate
provided that an accurate method was used to solve the original system.

The assumption of the initial solution being accurate enough breaks down
in many practical settings, for example when using threshold pivoting for
sparse or dense matrix factorizations~\citep{lindquist2022thrshpiv} or
even more complex error compensation
methods~\citep{lindquist2023beamaddmod}. In such situations, the solver
is often not as backward-stable as possible because it trades off its
performance for lower accuracy, which then has to be improved in the
iterative refinement stage.  Thus this second stage becomes the integral
part of the solution process and relies heavily on the quality of the
solution vector $x$ computed in the first stage. Another
aspect that we dive into shortly is the conditioning of the system
matrix $A$, which may be too high to allow a
less accurate solver resulting in the
iterative refinement being unable to
recover the desired accuracy. This is often
manifested by the iterations diverging or permanently stagnating.

Iterative refinement with artificially lowered precision was originally
advocated for the LU factorization~\citep{langou2006mxp}.
This was an example of purposefully
lowering the accuracy of the L and U factors of the LU factorization in
order to gain performance in the time-consuming portion of computing the
solution. To keep the accuracy in line with the expectations for the
high-precision solver, the iterative refinement was used to recover the
digits lost due to the use of a lower-precision. This is the
aforementioned change in the purpose of the solution refining
iterations: they became the method to bring the accuracy back to the
acceptable level rather than simply add a few digits of extra precision
beyond what is possible when using a uniform precision algorithm. The
motivating reason for this increased solver complexity was the
achievable performance gain: at the time it was asymptotically
$2\times$ faster for most commodity and high-end hardware and could
be even higher on
specialized hardware or in situation when a combination of hardware
arithmetic and software emulation were used. The accuracy could be
guaranteed with low values of the condition number of the original
system matrix $\kappa(A)$ relative to the unit roundoff of the lower
precision $u_{\text{low}}$ rather than the traditional requirement that
relies to conditioning with respect to the higher precision roundoff
$u_{\text{high}}$.

Introduction of GMRES~\citep{saad1986gmres} to the iterative refinement
class of methods lead to new benefits that could be gained both
numerically and in terms of performance.
Named GMRES-IR~\citep{carson2017iterref},
 this method added Krylov subspace iteration
thus replacing the original stationary iteration and doing away with its
main limitation: high sensitivity to the matrix conditioning. In fact,
ill-conditioning was the original design
constraint of GMRES-IR thus the
use of mixed-precision could be more expansive. This allowed much better
match with the hardware that was progressively fitted with a wider
variety of floating-point formats offering a large range of peak compute
levels. In fact, it was possible to use 3
precisions~\citep{carson2018iterref} or even
5 precisions~\citep{amestoy2024iterref5}.

In addition to using the right algorithmic framework to achieve
convergence in iterative refinement and desired accuracy of the solution
vector, there are still multiple issues to be resolved and they relate
mainly to properly done number scaling to fit the lower-precision
formats~\citep{luszczek2017half} and also preserve relevant matrix
properties~\citep{higham2019hlfprcsqz}. Note that numerical solvers
differ in this respect to the number scaling techniques employed for ML
and AI models, which use some form of microscaling that was initially
proposed for the large scale large language model workloads~\citep{rouhani2020micrsacl}.

The main reason for using the combination of mixed-precision and
iterative refinement is to adapt to hardware platform's compute
capabilities thus obtaining superior performance and energy consumption.
These aspects were extensively studied~\citep{haidar2018mxpiterref,haidar2020mxpiterref} on NVIDIA's GPU
accelerators along with more nuanced aspects that were introduced into
the solvers when Krylov subspace methods were
integrated into the mixed-precision
iterative refinement. Among others, these new considerations include
accelerating
convergence~\citep{baker2003gmrescnvacc, baker2005gmrescnvacc},
numerical stability~\citep{drkosova1995gmresnumstab,
paigeResidualBackwardError2001}, synchronization
cost~\citep{swirydowicz2019gmreslowsyn}, communication
pipelining~\citep{cools2019recurpipecg}, loss of
orthogonality~\citep{bjorck1992orthogrmsch,
paigeModifiedGramSchmidtMGS2006, paigeUsefulFormUnitary2009,
paigeEffectsLossOrthogonality2018}, and flexible variants for
non-constant
preconditioners~\citep{saadFlexibleInnerouterPreconditioned1993}.
Unfortunately, most of these topics are too broad to be within the scope
of this writing but might need to be explored further by the implementers
of the benchmark depending on the target hardware platform's processing
capabilities and data formats.

A successful launch and progress of a benchmark may be measured by the
scale and type of the machines,
on which it can be run. By that metric alone, HPL-MxP now was ported to
and subsequently optimized on the largest
supercomputers recorded on TOP500 list. In fact, it was one of the
very first codes that recorded Exascale-level performance with
$1$~Exa-OP recorded on a CPU-only Supercomputer
Fugaku~\citep{kudo2020exaflopfugaku} with only slightly over
$0.4$~Eflop/s recorded for HPL. A different implementation was required
to achieve scalable code for the largest machines with GPU
accelerators, and there are usually GPU-specific considerations to be
taken into account~\citep{lu2022summitfrontiermxp}. Such efforts led to
nearly $10$~Exa-OP of achieved performance while still maintaining the
accuracy of the HPL code. It is necessary to stress the difference
between the performance metrics of HPL and HPL-MxP since they use
drastically different floating-point formats that exhibit a range of
behaviors that may not be suitable for all users. However, the formula
to compute these figures of merit is identical for both codes:
$\frac{2}{3}n^3 \times t^{-1}$ with $n$ being the input matrix or
vectors sizes and $t$
being the time for solver phases including the iterative refinement in
case of HPL-MxP. At the same time, the precision and accuracy for both
the inputs and outputs is treated in the same manner thus giving a
similar numerical guarantees. Only the internal details of both compute
and storage for HPL-MxP remain the notable difference.



\section{Related Work}
\label{sec:related}

Mixed-precision hardware brings about not just new data formats but may
also violate the typical behavior of floating-point arithmetic.
For example, monotonic summation satisfies the property that for any
$x<\bar{x}$ and $y<\bar{y}$ implies
$\mbox{FP}(x+y)\le\mbox{FP}(\bar{x}+\bar{y})$, which is requested by the
IEEE 754 standard but may be violated in the recent
hardware~\citep{fasi2021nvtcore}.
With only mild assumptions about the monotonicity of the floating-point
arithmetic the bisection eigenvalue methods are accurate and numerically
stable~\citep{95ETNAbisecion}. This allows one
to recursively compute the
selected eigenvalues of a matrix with the \emph{bisection eigenvalue
algorithm}~\citep{95ETNAbisecion} that also applies for bidiagonal
matrices resulting from SVD~\citep{94DCSVD}.

When computing eigenvalues for some matrices with the bisection method,
the intermediate values produced in the Sturm function can grow or
shrink \emph{monotonically} as they are processed along the diagonal and
for the off-diagonal elements. Thus related to the mixed-precision
context, for the very large matrix sizes this could result in an
underflow or overflow as established by
mixed-precision eigensolvers~\citep{luszczek2024eigslice}.
Consequently, the remedial scaling techniques may be borrowed to fit the
limited exponent and mantissa of the lower-precision floating-point
representations.

Iterative refinement was originally suggested and analyzed
in numerical roundoff errors~\citep{wilkinson1963rea, wilkinson1965aep}
and was deemed beneficial
provided the accumulation of the residual vector entries are performed
in higher precisions. Then it was proposed~\citep{langou2006mxp}
instead to lower
the precision of the input matrix and to use the refinement
iteration for bringing back the lost digits. The potential
conditioning issues were subsequently addressed by introduction of
Krylov subspace methods~\citep{carson2017iterref} and three-precision
refinement as proposed earlier \citep{carson2018iterref}. This concept was
later taken further in direct sparse solvers~\citep{amestoy2024iterref5}
with the use of up to five precisions in refinement.

In a similar fashion, Cholesky QR is an example of an established method
that could leverage a mixed-precision approach. This allows utilization
of a higher precision format to maintain accuracy otherwise not possible
in a uniform precision scenario as observed
in numerical experiments~\citep{yamazaki2015cholqr}.

Within the field of Krylov subspace methods, GMRES received attention
early~\citep{hoggFastRobustMixedprecision2010}
and improved on newer hardware~\citep{lindquist2020mixprecgmres} so that
the mixed-precision could
benefit the time-to-solution but with a potential increase in the
iteration count due to stagnation of convergence when one of the
precisions was lowered excessively or the prerequisite properties of the
matrix were not met.

Mixed precision is also possible for finding eigenvalues and
eigenvectors of symmetric and Hermitian matrices as
shown~\citep{tsai2022eigmxp} for a variety of GPU accelerators with speedups
depending on the disparity between lower- and higher-precision hardware
units.

Mixed-precision was also proposed for Jacobi
SVD~\citep{ogita2020svditerrefmtx, uchino2024svditerref}
improving significantly upon the earlier methods.

Newton iteration can also be performed inside a mixed-precision
framework to achieve orthogonalization of eigenvector when the
lower-precision scheme failed to deliver sufficient numerical
quality~\citep{luszczek2024eigslice}.

A more extreme case of utilizing mixed-precision schemes is to switch
exclusively to integer compute units~\citep{ootomo2024intdgemm}. This is
performed with emulation of floating-point representation and can only
deliver performance benefits if the floating-point units are severely
under-provisioned in terms of compute capability.

Benchmarking of mixed-precision methods for sparse matrices was
proposed~\citep{yamazaki2022hpgmres} as an HPG-MxP benchmark or HPC GMRES in
mixed-precision. The reported results indicate speedups at scale on
large supercomputing installations.



\section{Algorithmic Framework}
\label{sec:algo}

The main purpose of running the HPL-MxP benchmark is to solve a system
of linear equations represented by a matrix $A$ to 64-bit floating-point
accuracy by first performing a mixed-precision factorization of a matrix
with the minimum precision required for the factors to be FP16 following
the available numerical analysis~\citep{carson2017iterref,carson2018iterref}.
Then an approximate solution $x_0$ is computed from the low-precision
factors $L$ and $U$. The standard LU decomposition can be used for this
task and no partial pivoting is required (but may be beneficial in a more
specific application setting) due to the way that the system matrix $A$
is generated as described below in the section on matrix generation.
Second, the approximate solution $x_0$ is used as the initial vector for
a non-stationary iterative method in 64-bit precision to generate
subsequent improvements to the initial low-precision approximation to
generate a final solution $x_{\text{final}}$ that should attain the
accuracy achieved by the LU decomposition in 64-bit floating-point
arithmetic such as the one that is computed by HPL.
GMRES~\citep{saad1986gmres,saadFlexibleInnerouterPreconditioned1993} is
the method recommended for the refinement's non-stationary iteration
because it is numerically
stable~\citep{carson2017iterref,carson2018iterref} when using the
low-precision LU factors as a left-side preconditioner in this iterative
algorithm.

The benchmark implementations should use the equivalent of the original
HPL benchmark harness~\citep{Dongarra2003} with a suitable modifications
in the matrix generator code, which needs to produce a non-symmetric
random matrix that admits non-singular LU decomposition in 64-bit
floating-point without partial pivoting. We explore this in more detail
in sections on matrix generation and experimental results. Note that the
matrices that produce diagonal elements slightly below the actual pivots
from partial pivoting would result in faster element growth and thus
require more iterations in the subsequent
refinement~\citep{liMakingSparseGaussian1998}.

Due to the number of built-in options to choose from when implementing
HPL-MxP there is an inherent need to require a consistent way of
obtaining reliably a figure of merit representing the result of running
a particular implementation of the benchmark. In order to obtain this
kind of uniformity across all computing platforms and make the reported
performance we require that the implemented algorithm for solving the
low-precision representation of the input system of equations must
numerically conform to an LU factorization. Specifically, we require
that the total operation count in that phase of the algorithm must
be on the order of $\frac{2}{3}n^3 + O(n^2)$ in terms of floating-point
operations even though 64-bit precision arithmetic is not required in
contrast to the rules for HPL.

Borrowing from the HPL's testing harness, HPL-MxP computes a
backward-error upon the completion of the solver and after the
convergence of the refinement iterations. The formula used for this
purpose is
$\|Ax-b\|_{\infty} \times (\|A\|_{\infty} \|x\|_{\infty} +
\|b\|_{\infty})^{-1} \times (n \times \epsilon)^{-1}$, where $\epsilon$
is the machine precision in 64-bit floating point arithmetic. On the
IEEE-compliant hardware $\epsilon$ is
$2^{-53}$ in binary~\citep{ieee1985fp754,ieee2019fp754}.
The use of $n$ or the size of the input matrix and vectors allows us to
remove the restriction on the system size. The correct solution is
assumed if the backward-error formula is less than $16$. Otherwise, the
run is deemed invalid. However, this form of the quality criterion and
the upper bound on the backward error are made to strictly conform with
the standard established by HPL so that the new mixed-precision method
with iterative refinement is guaranteed the same level of accuracy. For
a specific application, this may be relaxed or tightened depending on
the upstream requirements imposed on the linear solver.

As was already mentioned before, the lower-precision formats are much
more constrained than the 64-bit floating-point numbers.  To compensate,
the implementations are permitted to perform balancing and scaling of
the input matrix. The former is again borrowed directly from the HPL
testing harness and requires
$\Theta(n^2)$ operations to improve the
distribution of matrix elements' values. This directly affects the quality
of the L and U factors as well as the speed of convergence of the
refinement process.
The latter is more specific to mixed-precision
methods and is meant to transform the numbers so that they fit within
the range of the lower-precision floating-point formats chosen by the
implementation. Note however that the time consumed by both balancing
and scaling must be included in the total time to solution.

\begin{figure}
\centering
\begin{tabular}[c]{ | c | c | }
\hline
$A_{1,1}$ & \\
& \mbox{32 bits} \\
\hline
  & \\
32  & \\
  & \mbox{16 bits} \\
bits  & \\
  & \\
\hline
\end{tabular}
\caption{Computation of the Schur complement with higher-precision
devoted to accumulation of the trailing matrix.}
\label{fig:acc32}
\end{figure}

The implementation of the LU factorization is permitted to use a custom
mixed-precision scheme during the construction of the L and U factors.
In Figure~\ref{fig:acc32}, we show an example of such a scheme whereby
the panel factorization and triangular solves are done in 32-bit
arithmetic while the Schur complement or the update of the trailing
matrix is computed exclusively in 16-bit arithmetic but with
32-bit accumulation of the result.
This is performed with a matrix-matrix multiplication of the form
$A_{2,2} \leftarrow A_{2,2} - A_{2,1} \times A_{1,1}^{-1} \times
A_{1,2}$. This is possible on NVIDIA hardware that has
compute units capable of performing such an operation
at much higher rate than is the case for the larger data formats.

The figure of merit is the rate of execution computed exactly as the
number of canonical operations performed by HPL or $\frac{2}{3} n^3 +
\frac{3}{2} n^2$. This formula is composed of $\frac{2}{3} n^3 -
\frac{1}{2} n^2$ LU factorization operations and $2 n^2$ subsequent
operations for both back- and forward-substitution solves. The canonical
operation count is divided by the total time-to-solution, which then is
reported as the final rate of operations per second. They are not
the traditional flop/s since the floating-point operations
are likely to be different at different stages of the
computation.

Only the canonical operation count is used in the asymptotic compute
rate. However, the time-to-solution value accumulates many operations
required to compute an accurate solution vector. This could potentially
include a number of subsequent phases. Namely, the optional balancing
and/or scaling operations are performed first on the input matrix in
case overflow or underflow is suspected to occur, and it must be
prevented to improve refinement iteration convergence. Then the LU
factorization is computed in lower precision, which often would be
the most consuming part of the process since it would need the number of
operations roughly estimated by the canonical operations from the
formula for the asymptotic execution rate. Then the GMRES iterations
refine the solution to meet the accuracy criterion designed for the
64-bit floating-point arithmetic. This phase would include the use of
the lower-precision LU factors as a left-preconditioner. Only a handful
of iterations is expected to lead to convergence but the global
synchronizations could make this phase relatively long in duration. If
the particular implementation needs more than 50 iterations to converge
then the code should trigger a failure, and the run is considered
invalid. This tends to strike the balance between the quality of the
low-precision L and U factors and the performance of the refinement
process: only a handful of iterations are commonly needed in practice.
Finally, if either balancing or scaling was performed at the
beginning, then their effects have to be undone in order to project
it back onto the space of the input matrix. The time to compute the
backward error is not included and thus may be performed by an
unoptimized piece of code.

Conceptually, the solver is simple and results from a well researched
set of numerical methods. However, there are many crucial details that
may be lost by the implementer when considering the inherent complexity
of the modern supercomputing hardware platforms. For this reason we have
provided a reference implementation whose sole purpose is to show how
the benchmark code could look like initially before the
hardware-specific optimizations are attempted. By no means do we expect
this reference code to deliver near-optimal performance results, but at
the same time it should be possible to use it for estimation of a good
initial indicator of the hardware's capability. Further optimizations
should certainly be done to obtain a publishable result. As was the case
with other benchmarks~\citep{Dongarra2003, Luszczek2006, Dongarra2013,
Dongarra:2018:hpcg}, the hardware vendors and community provided such
optimized implementations that now cover CPUs and all major GPUs used in
the high end HPC systems.



\section{High Level Language Implementation}

There are increasing benefits of using higher-level languages to
implement mixed-precision solvers.
On the one hand, low-level programming languages such as the C
language~\citep{kernighan1978c} are
missing convenient support for \emph{generic programming}
with only a rudimentary functionality available by using C11 standard's
\textsf{\_Generic()} macro feature.
On the other hand, the C\kern-.17em\texttt{++}
language~\citep{ellis1990cplusplus}
offers higher level facilities to ease the generic interface design and
provides implementation mechanism through templates for multiple
precisions co-existing in a single code base while the
C\kern-.17em\texttt{++} concepts allow the implementers to provide
robust checking of the template parameters supplied by the users and
thus clearly indicate potential syntax errors concisely and consistently
as long as one of the newer language standards are used starting with
C\kern-.17em\texttt{++}20 onward.  In an even newer efforts, the
\textless\textsf{linalg}\textgreater~header file was added to the
standard library and was recently approved for inclusion in the ISO
C\kern-.17em\texttt{++}26 standard edition. The full specification as
proposed through the multiple subcommittee working versions of the
proposal are available as P1673~\citep{hoemmen2023linalg}. The
plethora of \emph{customization points} allow an easy adaptation of
these upcoming interface designs to our needs in a mixed-precision
iterative refinement setting.




\section{Matrix Generation}
\label{sec:mtxgen}

Scalable input data generation is essential for benchmarking large
supercomputing installations as was
observed repeatedly~\citep{Dongarra2003,
Luszczek2006, Dongarra2013, Dongarra:2018:hpcg,
luszczek2020scaldatagen},
and it was indicated that the lack of customizable inputs lead
to eventual stagnation of deployment leading
to irrelevance on constantly growing high-end systems. We will introduce
here a few potential solutions to this problem for HPL-MxP and later
provide experimental results supporting yet another matrix generation
scheme suitable for arbitrarily large benchmark runs. The main property
required from the input matrix is that, by construction, it should admit
non-singular L and U factors in 64-bit precision without the use of
pivoting for numerical stability. Matrices causing a large pivot growth
are clearly out of scope~\citep{turing1948rounderrmtx,zielke1974maxkond}.
Another requirement is to maintain spectral properties of the matrix
that makes it challenging for an unpreconditioned GMRES so that
generating a non-trivial preconditioner in the form of LU factorization
is a must. This is addressed by the generators described below.

Non-symmetric positive-definite matrices may be easily proven to not
require pivoting in the LU factorization. Examples of such matrices can
be easily generated from random matrices by scaling the diagonal
by the matrix size $n$. However, we consider such input matrices to
only serve as a good starting point for debugging and tuning purposes
but \emph{should not} be used in production runs for official results
because of how it distorts the matrix spectrum, causes underflow for
off-diagonal elements, produces easily predictable solution vectors
without the need of
performing many of the essential computations that form the basis of the
mixed-precision solver that requires non-stationary iterative refinement.

Instead of scaling the diagonal by the matrix size, a much more
constrained diagonal scaling was named \emph{departure from diagonal
dominance} and was proposed for a matrix
generator~\citep{luszczek2020scaldatagen}. This made the
generated matrices challenging for unpreconditioned GMRES so that no
spuriously fast convergence occurs in less than 50 iterations. In fact,
over 700 iterations were observed for larger matrices to achieve full
convergence without a preconditioner.

An alternative method was proposed~\citep{fasi2021tuinftycond} based
on the one-parameter matrix introduced
earlier~\citep{ostrowski1954oneparam}.
While the original definition was predictable, the further modifications
introduced randomization and controlled condition number thus regulating
suitability for the benchmark based on additional parameters. The only
drawback could be the necessity for a non-linear solver to arrive at
some of the matrix generator parameters.

Another option is to judiciously select either uniform or Gaussian
random distributions for the matrix elements because those can
positively regulate the spectrum of the generated matrices thus directly
affecting the convergence profile of a Krylov
solver~\citep{luszczek2020scaldatagen}.



\begin{figure*}[tb]
\centering
\includegraphics[width=.85\textwidth]{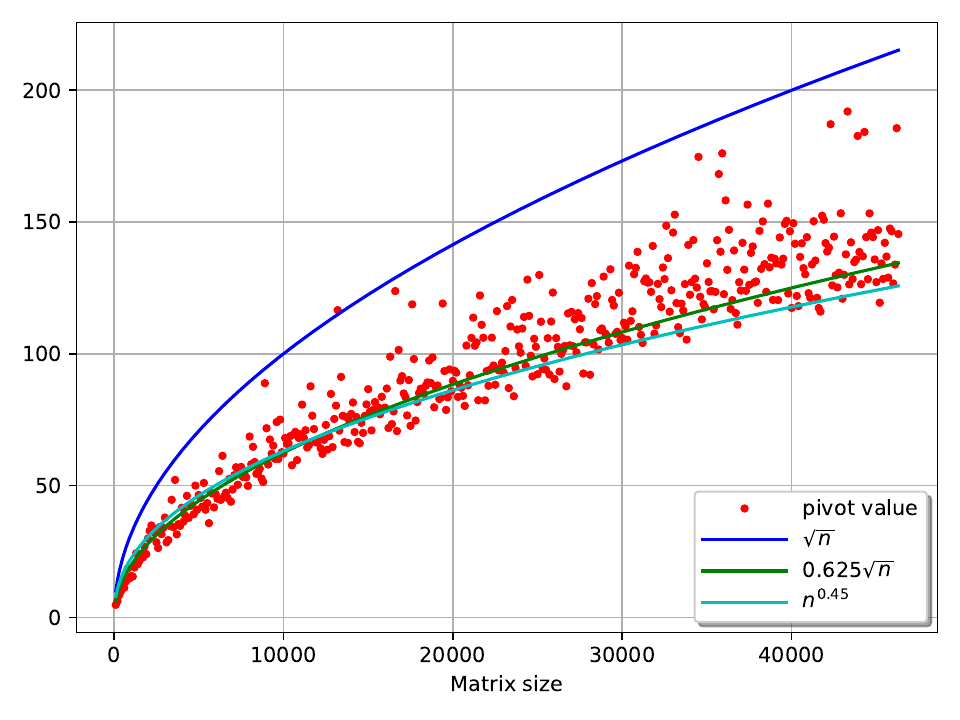}
\caption{The largest pivot value selected (red dots) for random
matrices of different sizes and their relative size compared to
theoretical growth curves.}
\label{fig:piv}
\end{figure*}

\section{Experimental Results}
\label{sec:results}

In this section, we show several of numerical experiments that attempt
to illustrate some of the numerical aspects of the computation of the LU
factorization that is essential for HPL-MxP and was not yet covered by
the references provided so far. We focus on the reproducibility aspects
in our experiments and thus we chose the input data sizes to easily fit
on most CPU and GPU systems while at the same time are large enough to
show relevant trends in numerical roundoff giving the reader an
opportunity to compare with the theoretical results provided in the
included references. This limited the sizes we used to $46340$ or the
largest value that fits in a 32-bit integer, and even 32-bit BLAS can be
used to produce the presented results.
The experimental results focus on matrices with uniformly
distributed random elements because the normally distributed entries are
better behaved in practice \citep{luszczek2020scaldatagen}.

As we mentioned earlier, in order to achieve an acceptable level of
numerical accuracy and pass the standard HPL criterion for results'
validity, it is necessary to first generate the input matrix that admits
non-singular L and U factors in 64-bit arithmetic without the use of
partial pivoting. At the same time we aim at only a mild pivot growth so
that the computed factors do not lose additional digits of accuracy,
which then would have to be recovered by the GMRES-based iterative
refinement. Finally, we also need to avoid excessive scaling of the
diagonal elements after generating uniformly random matrix elements so
that the diagonal and off-diagonal entries contribute in a balanced
fashion to the entries of the solution vector. To measure this trade
off, we measured the largest pivot value chosen by the standard LU
factorization in 64-bit floating-point arithmetic.
Figure~\ref{fig:piv} shows with red dots the pivots for a wide range of
matrix sizes for uniformly random matrices with values between $-1$ and
$1$ without scaling. In order to estimate the relevant pivot size, we also
included in the figure the analytical curves as a function of input
matrix size~$n$. The effective upper bound on the pivot value is
$\sqrt{n}$, which is in line with the random matrix
theory~\citep{edelman1989rndcond,higham2019problrndoff}.
However, the majority of the pivot
values are lower than the theoretical results suggest. In order to
estimate the relative size of these smaller pivots, we also show
analytical curves for $\frac{5}{8}\sqrt{n}$ and $n^{0.45}$. Both of them
are a good estimate for these smallest pivots, but the square root
function follows the trend better and should be preferred in the pivot
and scaling formulas.

It is also worth noting that the largest pivot value occurs in the last
$1\%$ of the columns, which corresponds well with the intuition that
the accumulation of trailing matrix updates would be the largest around the
very bottom end of the diagonal. We observed similar effects for
Wilkinson matrices in a completely different numerical
context~\citep{luszczek2024eigslice}.

\begin{figure*}[tb]
\centering
\includegraphics[width=.85\textwidth]{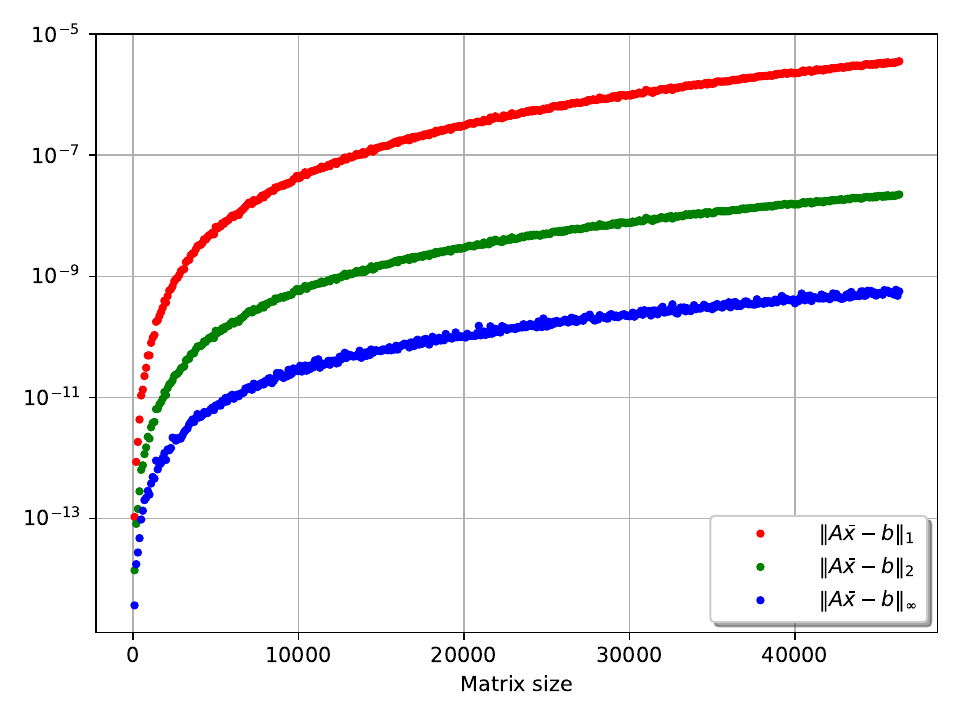}
\caption{Three norms of the solution vector of a solver with pivoting
for different matrix sizes.}
\label{fig:normpiv}
\end{figure*}

Figure \ref{fig:normpiv} shows the three common residual norms that could be
used for computing backward-error of the LU factorization with partial
pivoting in 64-bit floating-point arithmetic. The input matrix elements
were drawn from a uniformly random distribution of values between $-1$
and $1$ without any additional scaling of diagonal elements.
Also, these norms directly
determine if the criterion was satisfied to render the solution vector
as valid. Obviously, by their definitions, these norms are connected by
the scaling factors: $\sqrt{n}\|\cdot\|_{\infty} \le \|\cdot\|_2 \le
\sqrt{n^{-1}}\|\cdot\|_1$, which we can confirm with the figure's plots.
However, this is not the purpose of the figure. Instead, we point out
the variance values around the hypothetical trend line for the
computed norms across the tested matrix sizes. Even though all three
variance values may seem relatively small in the figure, we note
that the logarithmic scale is used for the $y$-axis, and thus the
visually small variance is still large in absolute terms. Also, the
biggest variance is associated with the smallest norm, which may cause
larger changes in the backward errors computed with that norm and thus
may potentially contribute to the failure of the
the passing criterion and render the
result invalid simply due to a natural variability in the input matrix
data. Finally, the figure will serve as the basis for comparing with the
LU factorization without any pivoting to give us a better idea about the
need for optimal scaling of the diagonal entries of the input matrix in
order to avoid an excessive element growth in the computed L and U
factors.

\begin{figure*}[tb]
\centering
\includegraphics[width=.85\textwidth]{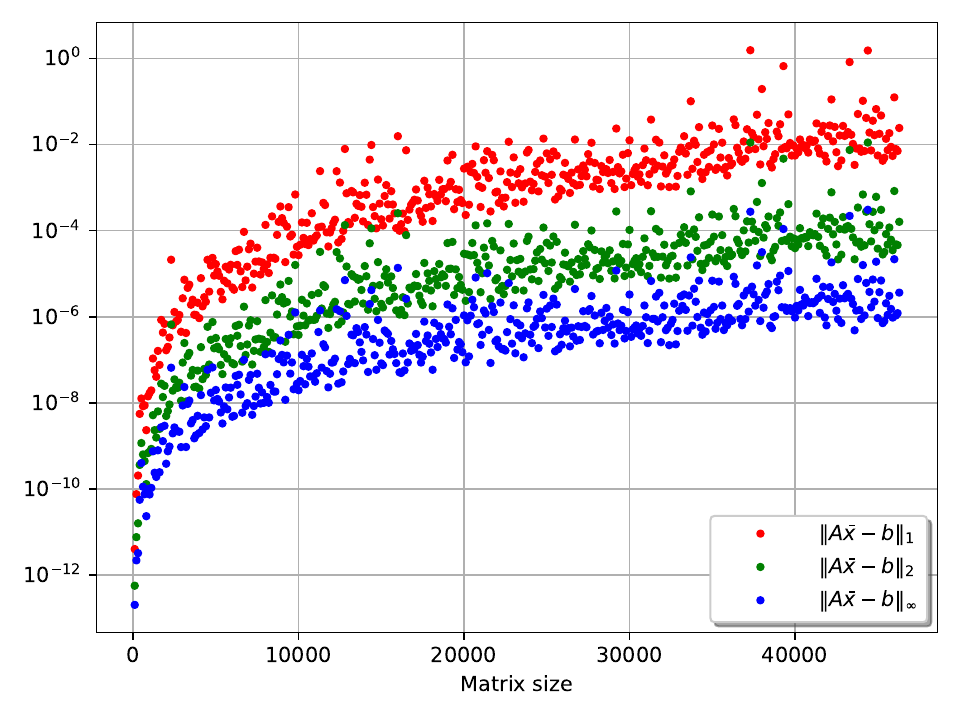}
\caption{Three norms of the solution vector of a solver without pivoting
for different matrix sizes.}
\label{fig:normnpv}
\end{figure*}

Figure \ref{fig:normnpv} shows the three common residual norms that could be
used for computing backward-error of the LU factorization \emph{without}
partial pivoting in 64-bit floating-point arithmetic. The input matrix
elements were again drawn from a uniformly random distribution of values
between $-1$ and $1$ without any additional scaling of diagonal
elements, which matches the input matrices from Figures~\ref{fig:piv}
and~\ref{fig:normpiv}. The combination of the lack of diagonal scaling
and not using partial pivoting contributes to catastrophic deterioration
of the values of the residual norms. We can estimate that each of the
norms increased by $2$ to $3$ orders of magnitude when minimum norm
values are compared, which constitutes an optimistic scenario and is not
that common for the tested input matrix sizes. The pessimistic cases have
norms larger by as much as $5$ orders of magnitude. And the
number of these worst-case occurrences is large enough to warrant
application of preventive measures. What contributes to such a
disparity between the best- and worst-case scenarios is greatly
increased variance of the norm values around the trend line.

The combined results presented in Figures~\ref{fig:piv}, \ref{fig:normpiv},
and~\ref{fig:normnpv} indicate that suitable scaling of the diagonal
elements of the input matrices is essential to achieve sufficiently
accurate L and U factors in 64-bit precision. The benefits would be
further compounded for even lower compute precisions. The conditioning
and roundoff error analysis based on random matrix theory suggest
the relevant scaling of $\sqrt{n}$, but it could be lowered at the
expense of more iterations in the mixed precision refinement process.



\section{Conclusions and Future Work}
\label{sec:conclude}

We presented the HPL-MxP benchmark with a backdrop of a combination of
mixed-precision and iterative refinement methods that have recently
received much attention from the numerical linear algebra research
community. We showed how these algorithms were proposed based on the
initial stationary iterative refinement design by lowering the precision
of intermediate computations and introducing Krylov subspace iteration.
We also presented a number of suitable choices for the scalable matrix
generator that makes partial pivoting in the LU factorization optional
while still allowing the L and U factors to serve as a
left-preconditioner for the GMRES iteration, which is a preferred way to
refine the initial solution to achieve full 64-bit accuracy irregardless
of the intermediate storage and compute formats. The reported results
have already achieved several Exa-operations per second and many
optimized implementations exist to allow the users to test a variety of
their systems based on CPUs, GPUs, or a combination of both.

One of the future research directions is to investigate the recent
advances in Generalized Random Butterfly Transform proposed
for the Gaussian Elimination~\citep{lindquist2024rbt}.
This may offer yet another scalable method
for generating input matrices that by construction do not require
pivoting.

Finally, we envision an ongoing effort of continuously improving the
reference implementation and collecting the results from the newly
introduced supercomputing systems. Improving the validation procedures
is also an important aspect for a benchmark that strives for maintaining
relevance and longevity.


\bibliographystyle{unsrt}
\bibliography{numhwsw}

\begin{thebibliography}{10}

\bibitem{moore1965cramming}
G.~E. Moore.
\newblock Cramming more components onto integrated circuits.
\newblock {\em Electronics}, 38(8), 1965.

\bibitem{Dongarra2003}
Jack~J. Dongarra, Piotr Luszczek, and Antoine Petitet.
\newblock The {LINPACK} benchmark: Past, present, and future.
\newblock {\em Concurrency and Computation: Practice and Experience},
  15(9):803--820, August 10 2003.
\newblock DOI: 10.1002/cpe.728. ISSN 1532-0634.

\bibitem{carson2018iterref}
Erin Carson and Nicholas~J. Higham.
\newblock Accelerating the solution of linear systems by iterative refinement
  in three precisions.
\newblock {\em SIAM J. Sci. Comput.}, 40(2):A817--A847, 2018.

\bibitem{carson2017iterref}
Erin Carson and Nicholas~J. Higham.
\newblock A new analysis of iterative refinement and its application to
  accurate solution of ill-conditioned sparse linear systems.
\newblock {\em SIAM J. Sci. Comput.}, 39(6):A2834–A2856, 2017.

\bibitem{haidar2018mxpiterref}
Azzam Haidar, Stanimire Tomov, Jack Dongarra, and Nicholas~J. Higham.
\newblock Harnessing {GPU} tensor cores for fast {FP16} arithmetic to speed up
  mixed-precision iterative refinement solvers.
\newblock In {\em Proceedings of the International Conference for High
  Performance Computing, Networking, Storage, and Analysis}, SC '18. IEEE
  Press, 2018.

\bibitem{haidar2020mxpiterref}
A.~Haidar, H.~Bayraktar, S.~Tomov, J.~Dongarra, and N.~J. Higham.
\newblock Mixed-precision iterative refinement using tensor cores on gpus to
  accelerate solution of linear systems.
\newblock {\em Proceedings of the Royal Society A}, 476(2243), November 2020.

\bibitem{wilkinson1963rea}
James~H. Wilkinson.
\newblock {\em Rounding Errors in Algebraic Processes}.
\newblock Notes on Applied Science No. 32, Her Majesty's Stationery Office,
  London, 1963.
\newblock Also published by Prentice-Hall, Englewood Cliffs, NJ, USA. Reprinted
  by Dover, New York, 1994.

\bibitem{wilkinson1965aep}
James~H. Wilkinson.
\newblock {\em The Algebraic Eigenvalue Problem}.
\newblock pub-OUP, 1965.

\bibitem{lindquist2022thrshpiv}
Neil Lindquist, Mark Gates, Piotr Luszczek, and Jack Dongarra.
\newblock Threshold pivoting for dense {LU} factorization.
\newblock In {\em ScalAH22: 13th Workshop on Latest Advances in Scalable
  Algorithms for Large-Scale Heterogeneous Systems}, pages 1--10, Dallas,
  Texas, USA, November 14, 2022.

\bibitem{lindquist2023beamaddmod}
Neil Lindquist, Piotr Luszczek, and Jack Dongara.
\newblock Using additive modifications in {LU} factorization instead of
  pivoting.
\newblock In {\em International Conference on Supercomputing (ICS) 2023}, pages
  14--24, Orlando, Florida, USA, June 2023.
\newblock Best paper nominee.

\bibitem{langou2006mxp}
Julie Langou, Julien Langou, Piotr Luszczek, Jakub Kurzak, Alfredo Buttari, and
  Jack Dongarra.
\newblock Exploiting the performance of 32 bit floating point arithmetic in
  obtaining 64 bit accuracy (revisiting iterative refinement for linear
  systems).
\newblock In {\em Proceedings of the 2006 ACM/IEEE Conference on
  Supercomputing}, SC '06, page 113–es, New York, NY, USA, 2006. Association
  for Computing Machinery.

\bibitem{saad1986gmres}
Y.~Saad and M.~H. Schultz.
\newblock {GMRES:} a generalized minimal residual algorithm for solving
  nonsymmetric linear systems.
\newblock {\em SIAM J. Sci. Stat. Comput.}, 7:856–869, 1986.

\bibitem{amestoy2024iterref5}
Patrick~R. Amestoy, Alfredo Buttari, Nicholas~John Higham, Jean-Yves
  L'Excellent, Th{\'e}o Mary, and Bastien Vieubl{\'e}.
\newblock Five-precision gmres-based iterative refinement.
\newblock {\em SIAM J. Matrix Anal. Appl.}, 45:529--552, 2024.

\bibitem{luszczek2017half}
Piotr Luszczek, Jakub Kurzak, Ichitaro Yamazaki, and Jack Dongarra.
\newblock Towards numerical benchmark for half-precision floating point
  arithmetic.
\newblock In {\em 2017 IEEE High Performance Extreme Computing Conference
  (HPEC)}, Waltham, Massachusetts, USA, 2017.

\bibitem{higham2019hlfprcsqz}
Nicholas~J. Higham, Srikara Pranesh, and Mawussi Zounon.
\newblock Squeezing a matrix into half precision, with an application to
  solving linear systems.
\newblock {\em SIAM Journal on Scientific Computing}, 41(4):A2536--A2551, 2019.

\bibitem{rouhani2020micrsacl}
Bita Rouhani, Daniel Lo, Ritchie Zhao, Ming Liu, Jeremy Fowers, Kalin
  Ovtcharov, Anna Vinogradsky, Sarah Massengill, Lita Yang, Ray Bittner,
  Alessandro Forin, Haishan Zhu, Taesik Na, Prerak Patel, Shuai Che, Lok~Chand
  Koppaka, Xia Song, Subhojit Som, Kaustav Das, Saurabh Tiwary, Steve
  Reinhardt, Sitaram Lanka, Eric Chung, and Doug Burger.
\newblock Pushing the limits of narrow precision inferencing at cloud scale
  with microsoft floating point.
\newblock In {\em 34th Conference on Neural Information Processing Systems
  (NeurIPS 2020)}, Vancouver, Canada, 2020.

\bibitem{baker2003gmrescnvacc}
A.~H. Baker.
\newblock {\em On Improving the Performance of the Linear Solver restarted
  {GMRES}}.
\newblock PhD thesis, University of Colorado, 2003.

\bibitem{baker2005gmrescnvacc}
A.~H. Baker, E.R. Jessup, and T.~Manteuffel.
\newblock A technique for accelerating the convergence of restarted {GMRES}.
\newblock {\em SIAM J. Matrix Anal. Appl.}, 26(962), 2005.

\bibitem{drkosova1995gmresnumstab}
Jitka Drko\v{s}ov\v{a}, Anne Greenbaum, Miroslav Roslo\v{z}n\'ik, and
  Zden\v{e}k Strako\v{s}.
\newblock Numerical stability of {GMRES}.
\newblock {\em BIT Numerical Mathematics}, 35:309--330, September 1995.

\bibitem{paigeResidualBackwardError2001}
Christopher~C. Paige and Zdenvek Strakos.
\newblock Residual and backward error bounds in minimum residual {{Krylov}}
  subspace methods.
\newblock {\em SIAM J. Sci. Comput.}, 23(6):1898--1923, June 2001.

\bibitem{swirydowicz2019gmreslowsyn}
K.~\'Swirydowicz, J.~Langou, Shreyas Ananthan, U.~Yang, and S.~Thomas.
\newblock Low synchronization {Gram-Schmidt} and {GMRES} algorithms.
\newblock {\em Numerical Linear Algebra with Applications}, 2019.

\bibitem{cools2019recurpipecg}
Siegfried Cools, Jeffrey Cornelis, and Wim Vanroose.
\newblock Numerically stable recurrence relations for the communication hiding
  pipelined conjugate gradient method.
\newblock {\em IEEE Transactions on Parallel and Distributed Systems},
  30(11):2507--2522, 2019.

\bibitem{bjorck1992orthogrmsch}
{\r{A}}ke Bj\"{o}rck and Christopher~C. Paige.
\newblock Loss and recapture of orthogonality in the modified {Gram-Schmidt}
  algorithm.
\newblock {\em SIAM Journal on Matrix Analysis and Applications},
  13(1):176--190, 1992.

\bibitem{paigeModifiedGramSchmidtMGS2006}
Christopher~C. Paige, Miroslav. Rozlozn{\'i}k, and Zdenvek. Strakos.
\newblock Modified {Gram-Schmidt} {(MGS)}, least squares, and backward
  stability of {MGS-GMRES}.
\newblock {\em SIAM J. Matrix Anal. Appl.}, 28(1):264--284, 2006.

\bibitem{paigeUsefulFormUnitary2009}
Christopher~C. Paige.
\newblock A useful form of unitary matrix obtained from any sequence of unit
  2-norm n-vectors.
\newblock {\em SIAM J. Matrix Anal. Appl.}, 31(2):565--583, May 2009.

\bibitem{paigeEffectsLossOrthogonality2018}
Christopher~C. Paige.
\newblock The effects of loss of orthogonality on large scale numerical
  computations.
\newblock In Osvaldo Gervasi, Beniamino Murgante, Sanjay Misra, Elena Stankova,
  Carmelo~M. Torre, Ana Maria~A.C. Rocha, David Taniar, Bernady~O. Apduhan,
  Eufemia Tarantino, and Yeonseung Ryu, editors, {\em Computational Science and
  Its Applications -- {ICCSA} 2018}, pages 429--439, {Cham}, 2018. {Springer
  International Publishing}.

\bibitem{saadFlexibleInnerouterPreconditioned1993}
Youcef Saad.
\newblock A flexible inner-outer preconditioned {GMRES} algorithm.
\newblock {\em SIAM Journal on Scientific Computing}, 14(2):461--469, March
  1993.

\bibitem{kudo2020exaflopfugaku}
S.~Kudo, K.~Nitadori, T.~Ina, and T.~Imamura.
\newblock Implementation and numerical techniques for one {Eflop/s} {HPL-AI}
  benchmark on {Fugaku}.
\newblock In {\em ScalA20 Conf., SC20}, 2020.

\bibitem{lu2022summitfrontiermxp}
Hao Lu, Michael Matheson, Vladyslav Oles, Austin Ellis, Wayne Joubert, and
  Feiyi Wang.
\newblock Climbing the summit and pushing the frontier of mixed precision
  benchmarks at extreme scale.
\newblock In {\em SC22: International Conference for High Performance
  Computing, Networking, Storage and Analysis}, pages 1--15, 2022.

\bibitem{fasi2021nvtcore}
Massimiliano Fasi, Ncholas~J. Higham, S.~Pranesh, and M.~Mikaitis.
\newblock Numerical behavior of {NVIDIA} tensor cores.
\newblock {\em PeerJ Computer Science}, 7(e330), 2021.

\bibitem{95ETNAbisecion}
James~W. Demmel, Inderjit Dhillon, and Huan Ren.
\newblock On the correctness of some bisection-like parallel eigenvalue
  algorithms in floating point arithmetic.
\newblock {\em Electronic Transactions on Numerical Analysis}, 3:116--149,
  December 1995.

\bibitem{94DCSVD}
M.~Gu, J.~Demmel, and I.~Dhillon.
\newblock Efficient computation of the {Singular} {Value} {Decomposition} with
  applications to least squares problems.
\newblock Technical Report CS-94-257, Department of Computer Science,
  University of Tennessee, October 1994.

\bibitem{luszczek2024eigslice}
Piotr Luszczek, Anthony Castaldo, Yaohung~M Tsai, Daniel Mishler, and Jack
  Dongarra.
\newblock Numerical eigen-spectrum slicing, accurate orthogonal eigen-basis,
  and mixed-precision eigenvalue refinement using {OpenMP} data-dependent tasks
  and accelerator offload.
\newblock {\em The International Journal of High Performance Computing
  Applications}, 38(6):671--691, 2024.

\bibitem{yamazaki2015cholqr}
Ichitaro Yamazaki, Stanimire Tomov, and Jack Dongarra.
\newblock Mixed-precision {Cholesky} {QR} factorization and its case studies on
  multicore {CPU} with multiple {GPUs}.
\newblock {\em SIAM J. Sci. Comput.}, 37(3):C307–C330, 2015.

\bibitem{hoggFastRobustMixedprecision2010}
J.~D. Hogg and J.~A. Scott.
\newblock A fast and robust mixed-precision solver for the solution of sparse
  symmetric linear systems.
\newblock {\em ACM Trans. Math. Softw.}, 37(2):17:1--17:24, April 2010.

\bibitem{lindquist2020mixprecgmres}
Neil Lindquist, Piotr Luszczek, and Jack Dongarra.
\newblock Improving the performance of the {GMRES} method using mixed-precision
  techniques.
\newblock In {\em Proceedings of Smokey Mountain Conference, virtual}, 2020.

\bibitem{tsai2022eigmxp}
Yaohung Tsai, Piotr Luszczek, and Jack Dongarra.
\newblock Mixed-precision algorithm for finding selected eigenvalues and
  eigenvectors of symmetric and {Hermitian} matrices.
\newblock In {\em ScalAH22: 13th Workshop on Latest Advances in Scalable
  Algorithms for Large-Scale Heterogeneous Systems}, pages 1--10, Dallas,
  Texas, USA, November 14, 2022.

\bibitem{ogita2020svditerrefmtx}
T.~Ogita and K.~Aishima.
\newblock Iterative refinement for singular value decomposition based on matrix
  multiplication.
\newblock {\em J. Comput. Appl. Math}, 369:112512, 2020.

\bibitem{uchino2024svditerref}
Yuki Uchino, Takeshi Terao, and Katsuhisa Ozaki.
\newblock Acceleration of iterative refinement for singular value
  decomposition.
\newblock {\em Numerical Algorithms}, 95:979--1009, 2024.

\bibitem{ootomo2024intdgemm}
H.~Ootomo, Katsuhisa Ozaki, and Rio Yokota.
\newblock {DGEMM} on integer matrix multiplication unit.
\newblock {\em The International Journal of High Performance Computing
  Applications}, 2024.
\newblock in press.

\bibitem{yamazaki2022hpgmres}
Ichitaro Yamazaki, Christian Glusa, Jennifer Loe, Piotr Luszczek, Sivasankaran
  Rajamanickam, and Jack Dongarra.
\newblock High-performance {GMRES} multi-precision benchmark: Design,
  performance, and challenges.
\newblock In {\em 2022 International Workshop on Performance Modeling,
  Benchmarking and Simulation of High Performance Computer Systems (PMBS)},
  pages 1--10, Dallas, Texas, USA, November 14, 2022.

\bibitem{liMakingSparseGaussian1998}
Xiaoye~S. Li and James~W. Demmel.
\newblock Making sparse {{Gaussian}} elimination scalable by static pivoting.
\newblock In {\em {SC} '98: Proceedings of the 1998 {ACM}/{IEEE} Conference on
  Supercomputing}, pages 34--34, {San Jose, CA, USA}, November 1998.

\bibitem{ieee1985fp754}
Standard for binary floating point arithmetic, 1985.

\bibitem{ieee2019fp754}
{IEEE} standard for floating-point arithmetic, 2019.
\newblock (Revision of IEEE 754-2008) DOI:10.1109/IEEESTD.2019.8766229.

\bibitem{Luszczek2006}
Piotr Luszczek, Jack Dongarra, and Jeremy Kepner.
\newblock Design and implementation of the {HPCC} benchmark suite.
\newblock {\em CT Watch Quarterly}, 2(4A):18--23, November 2006.

\bibitem{Dongarra2013}
Jack Dongarra and Piotr Luszczek.
\newblock {HPC} {Challenge:} design, history, and implementation highlights.
\newblock In Jefferey~S. Vetter, editor, {\em Contemporary High Performance
  Computing: From Petascale Toward Exascale}, volume~1 of {\em CRC
  Computational Science Series}, chapter~2, pages 13--32. Taylor and Francis,
  Boca Raton, 1 edition, 2013.

\bibitem{Dongarra:2018:hpcg}
Jack Dongarra, Michael~A. Heroux, and Piotr Luszczek.
\newblock The {High-Performance} {Conjugate} {Gradients} benchmark.
\newblock {\em SIAM News}, 51(1):12--12, January/February 2018.

\bibitem{kernighan1978c}
Brian~W. Kernighan and Dennis~M. Ritchie.
\newblock {\em The {C} Programming Language}.
\newblock Prentice-Hall, Upper Saddle River, New Jersey, 1978.

\bibitem{ellis1990cplusplus}
Margaret~A. Ellis and Bjarne Stroustrup.
\newblock {\em The Annotated {C\kern-.17em\texttt{++}} Reference Manual}.
\newblock Addison-Wesley, Reading, Massachusetts, 1990.

\bibitem{hoemmen2023linalg}
Mark Hoemmen, Dasy Hollman, Christian Trott, Daniel Sunderland, Nevin Liber,
  Alicia Klinvex, Li-Ta Lo, Damien Lebrun-Grandie, Graham Lopez, Peter Caday,
  Sarah Knepper, Piotr Luszczek, and Timothy Costa.
\newblock {P1673R12:} a free function linear algebra interface based on the
  {BLAS}.
\newblock {\em ISO JTC21/SC22/WG21 Library Evolution Working Group}, pages
  1--141, March 14, 2023.

\bibitem{luszczek2020scaldatagen}
Piotr Luszczek, Yaohung Tsai, Neil Lindquist, Hartwig Anzt, and Jack Dongarra.
\newblock Scalable data generation for evaluating mixed-precision.
\newblock In {\em Proceedings of IEEE HPEC'20: High Performance Extreme
  Computing}, Waltham, MA, USA, 2020.

\bibitem{turing1948rounderrmtx}
Alan~M. Turing.
\newblock Rounding-off errors in matrix processes.
\newblock {\em The Quarterly Journal of Mechanics and Applied Mathematics},
  1(1):287--308, January 1948.

\bibitem{zielke1974maxkond}
G.~Zielke.
\newblock {Testmatrizen mit maximaler Konditionszahl}.
\newblock {\em Computing}, 13(1):33--54, March 1974.

\bibitem{fasi2021tuinftycond}
Massimiliano Fasi and Nicholas~J. Higham.
\newblock Matrices with tunable infinity-norm condition number and no need for
  pivoting in lu factorization.
\newblock {\em SIAM Journal on Matrix Analysis and Applications},
  42(1):417--435, 2021.

\bibitem{ostrowski1954oneparam}
A.~M. Ostrowski.
\newblock On the spectrum of a one-parametric family of matrices.
\newblock {\em J. Reine Angew. Math.}, 193(3/4):143–160, 1954.

\bibitem{edelman1989rndcond}
Alan Edelman.
\newblock {\em Eigenvalues and Condition Numbers of Random Matrices}.
\newblock PhD thesis, Massachusetts Institute of Technology, May 1989.

\bibitem{higham2019problrndoff}
Ncholas~J. Higham and Theo Mary.
\newblock A new approach to probabilistic rounding error analysis.
\newblock {\em SIAM journal on Scientific Computing}, 41(5):A2815–A2835,
  2019.

\bibitem{lindquist2024rbt}
Neil Lindquist, Piotr Luszczek, and Jack Dongarra.
\newblock Generalizing {Random Butterfly Transforms} to arbitrary matrix sizes.
\newblock {\em ACM Trans. Math. Softw.}, 50(4), December 2024.

\end{thebibliography}

\end{document}